\documentclass[12pt,reqno]{amsart}
\usepackage{graphicx, amssymb, color,pdfsync}
\usepackage[all,cmtip]{xy}
\numberwithin{equation}{section}
\usepackage{mathrsfs}
\usepackage{hyperref}

\hypersetup{
    colorlinks=true,       % false: boxed links; true: colored links
    linkcolor=blue,          % color of internal links
    citecolor=blue,        % color of links to bibliography
    filecolor=blue,      % color of file links
    urlcolor=blue           % color of external links
}
\usepackage{tikz}
\usetikzlibrary{matrix,arrows}

\usepackage[switch]{lineno}
\usepackage{yfonts}
\advance\hoffset by -0.9 truecm

\begin{document}
\newcommand{\s}{\vspace{0.2cm}}

\newtheorem{theo}{Theorem}
\newtheorem{prop}{Proposition}
\newtheorem{coro}{Corollary}
\newtheorem{lemm}{Lemma}
\newtheorem{claim}{Claim}
\newtheorem{example}{Example}
\theoremstyle{remark}
\newtheorem*{rema}{\it Remarks}
\newtheorem*{rema1}{\it Remark}
\newtheorem*{defi}{\it Definition}
\newtheorem*{theo*}{\bf Theorem}
\newtheorem*{coro*}{Corollary}

\title[A family of $(p,n)$-gonal Riemann surfaces]{A family of $(p,n)$-gonal Riemann surfaces with \\several $(p,n)$-gonal groups}
\date{}

\author{Sebasti\'an Reyes-Carocca}
\address{Departamento de Matem\'atica y Estad\'istica, Universidad de La Frontera, Avenida Francisco Salazar 01145, Temuco, Chile.}
\email{sebastian.reyes@ufrontera.cl}

\thanks{Partially supported by Fondecyt Grants 11180024, 1190991 and Redes Grant 2017-170071}

\begin{abstract} Let $p \geqslant 3$ be a prime number and let $n \geqslant 0$ be  an integer such that $p-1$ divides $n.$ In this short note we construct a family of $(p,n)$-gonal  Riemann surfaces of maximal genus $2np+(p-1)^2$  with more than one $(p,n)$-gonal group.

\end{abstract}
\maketitle

\section{Introduction and statement of the result}

Let $S$ be a compact Riemann surface of genus $g \geqslant 2$ and let $\mbox{Aut}(S)$ denote its automorphism group. If $p \geqslant 2$ is a prime number and $n \geqslant 0$ is an integer then $S$ is called {\it $(p,n)$-gonal} if  there exists a group of automorphisms $$ C_p \cong H \leqslant \mbox{Aut}(S)$$such that the corresponding orbit space $S/H$ has genus $n.$ The group $H$ is called a {\it $(p,n)$-gonal group} of $S.$

\s

Each compact Riemann surface with non-trivial automorphisms is $(p,n)$-gonal for suitable of $p$ and $n.$ This simple fact shows that  to study $(p,n)$-gonal Riemann surfaces and their automorphisms is equivalent to study the singular locus of the moduli space of compact Riemann surfaces.

\s

$(p,n)$-gonal Riemann surfaces and their automorphisms  have been extensively considered over the last century as they generalize important and well-studied classes of Riemann surfaces, such as $(2,0)$-gonal or {\it hyperelliptic}, $(p,0)$-gonal or {\it $p$-gonal} and $(2,n)$-gonal or {\it $n$-hyperelliptic} Riemann surfaces, among others.

\s

Let $S$ be a $p$-gonal Riemann surface of genus $g \geqslant 2.$
By the classical Castelnuovo-Severi  inequality (see Accola's book \cite{accola}), if \begin{equation} \label{ec}g >  (p-1)^2\end{equation}then the $p$-gonal group is unique  in the automorphism group of $S$. A family of $p$-gonal Riemann surfaces of maximal genus $g=(p-1)^2$ endowed with two $p$-gonal groups was constructed in \cite{ciy}, showing that the bound \eqref{ec} is sharp. Furthermore, in the general case, following \cite{gabino}, if $S$ has two $p$-gonal groups then they are  conjugate in the automorphism group of $S$. An upper bound for the number of such groups was obtained in \cite{grego}. 
 
 \s

For $(p,n)$-gonality with $n \geqslant 1,$ the Castelnuovo-Severi inequality ensures that if $S$ is a $(p,n)$-gonal Riemann surface of genus $g\geqslant 2$ and \begin{equation} \label{ec2}g > 2pn + (p-1)^2\end{equation}then the $(p,n)$-gonal group is unique in the automorphism group of $S$.
In the general case, it was proved in \cite{greww} that if $S$ is a $(p,n)$-gonal Riemann surface of genus $g$ and $p > 2n + 1$ then all its $(p, n)$-gonal groups are conjugate in the automorphism group of $S;$ an upper bound for the size of the corresponding conjugacy class was also determined in the same paper. Later,  in \cite{andreas}, the uniqueness of the $(p,n)$-gonal group was proved to be true under the assumptions that the $(p,n)$-gonal group acts with fixed points and $p > 6n-6.$ 
\s

This short note is devoted to provide a family of $(p,n)$-gonal Riemann surfaces of maximal genus $g=2pn+(p-1)^2$  with two $(p,n)$-gonal groups. The existence of this family shows that the bound \eqref{ec2} is sharp, for infinitely many pairs $(p,n).$

\begin{theo*} \label{t1} Let $p \geqslant 3$ be a prime number and let $n \geqslant 0$ be  an integer such that $p-1$ divides $n.$ Set $$d=n/(p-1) + 1.$$Then there exists a complex $d$-dimensional family of $(p,n)$-gonal Riemann surfaces $S$ of genus $$g=2np+(p-1)^2$$with automorphism group of order $4p^2$ acting on $S$ with signature $$(0; 2,2,2,p, \stackrel{d}{\ldots},p)$$in such a way that each $S$ has more than one $(p,n)$-gonal group.

\end{theo*}

\begin{rema1}
It is worth mentioning here the following observations which will follow from the proof of the theorem.
\begin{enumerate}
\item The result remains true if $p=2$ and $n$ is odd. 
\item If $n=0$ our family agrees with the family constructed in \cite{ciy}.
\end{enumerate}
\end{rema1}

\section{Proof of the Theorem}

Let $\Delta$ be a Fuchsian group of signature $(0; 2,2,2,p, \stackrel{d}{\ldots},p)$  canonically presented as $$\Delta=\langle \gamma_1,  \ldots, \gamma_{d+3} : \gamma_1^2=\gamma_2^2=\gamma_3^2=\gamma_4^p=\cdots =\gamma_{d+3}^p=\gamma_1\cdots\gamma_{d+3}=1\rangle$$and consider the group $G=\mathbf{D}_p \times \mathbf{D}_p$ (where $\mathbf{D}_p$ denotes the dihedral group of order $2p$)  presented in terms of generators $s_1, s_2, r_1, r_2$ and relations$$s_1^2=s_2^2=r_1^p=r_2^p=(s_1r_1)^2=(s_2r_2)^2=[s_1, r_2]=[s_1, s_2]=[r_1, r_2]=[r_1, s_2]=1.$$

{\it Existence of the family.} By virtue of the classical Riemann's existence theorem, the existence of the desired family follows after verifying that the Riemann-Hurwitz formula holds and after providing a surface-kernel epimorphism $\theta$ from $\Delta$ onto $G.$

\s

Note that the equality $$2(g-1)=4p^2(-2+3(1-\tfrac{1}{2})+d(1-\tfrac{1}{p}))$$shows that the Riemann-Hurwitz formula is satisfied for a branched $4p^2$-fold  covering map from a compact Riemann surface of genus $g=2np+(p-1)^2$ onto the projective line, ramified over three values marked with 2 and $d$ values marked with $p$. 

\s
In addition, if $d$ is odd we can choose the surface-kernel epimorphism $\theta$ as \begin{displaymath}
\theta(\gamma_1)= s_1, \,\, \theta(\gamma_2)=s_2, \,\, \theta(\gamma_3)= s_1s_2r_1r_2 \, \mbox{ and }\,  \theta(\gamma_i)= \left\{ \begin{array}{ll}
 (r_1r_2)^{-1} & \textrm{if $i$ is even}\\
\,\,\, r_1r_2 & \textrm{if $i$ is odd,}
  \end{array} \right.
\end{displaymath}and if $d$ is even we can choose $\theta$ as \begin{displaymath}
\theta(\gamma_1)= s_1, \,\, \theta(\gamma_2)= s_2, \,\, \theta(\gamma_3)=s_1s_2(r_1r_2)^{-d/2} \, \mbox{ and }\, \theta(\gamma_i)= \left\{ \begin{array}{ll}
r_1 & \textrm{if $i$ is even}\\
 r_2 & \textrm{if $i$ is odd}
  \end{array} \right.
\end{displaymath}where $i \in \{4, \ldots, d+3\}.$ 
\s

The complex dimension of the family agrees with the complex dimension of the Teichm\"{u}ller space associated to $\Delta;$ namely, its dimension is $d.$
\s

{\it $(p,n)$-gonal groups.} We denote the branched regular covering map given by the action of $G$ on $S$  by $\pi: S \to S/G$ and  its branch values by  $y_1,y_2, y_3, z_1 \ldots, z_{d}$,  where each $y_k$ is marked with 2 and each $z_k$ is marked with $p.$

\s

Assume $d$ odd. Consider the cyclic subgroups of order $p$ $$H_1 = \langle r_1r_2 \rangle \,\, \mbox{ and }\,\, H_2 = \langle r_1^{-1}r_2 \rangle$$of $G.$ We denote by $\pi_1$ and $\pi_2$ the  branched regular covering maps given by the action of $H_1$ and $H_2$ on $S$ respectively. We observe that the fiber of $\pi$ over each $y_k$ does not contain any branch value of $\pi_1$ and $\pi_2.$ In addition,  for each $k$,  the fiber of $\pi$ over $z_k$ has $4p$ elements; the isotropy group of $2p$ of them is  isomorphic to $H_1$ and the remaining ones have isotropy group isomorphic $H_2.$ It follows that $\pi_1$ and $\pi_2$ ramify over $2pd$ values, each of them marked with $p.$ Equivalently,  the signature of the action of $H_j$ on $S$ is $(n_j; p, \stackrel{2dp}{\ldots}, p)$ where $n_j$ is the genus of $S/H_j$. We now consider the Riemann-Hurwitz formula to see that $$2(g-1)=p[2n_j-2+2pd(1-\tfrac{1}{p})]$$and, after straightforward computations, one obtains that $n_j=n$ for $j=1,2$. 

\s

Assume $d$ even. Consider the cyclic subgroups of order $p$ $$H_1 = \langle r_1 \rangle \,\, \mbox{ and }\,\, H_2 = \langle r_2 \rangle$$of $G$ and let $\pi_1$ and $\pi_2$ be as before. As in the previous case, the fiber of $\pi$ over each $y_k$ does not contain any branch value of $\pi_1$ and $\pi_2.$ For each $k$ the fiber of $\pi$ over $y_k$ has $4p$ elements; the isotropy group  of them is isomorphic to $H_1$ if $k$ is odd and  is  isomorphic to $H_2$ if $k$ is even. It follows that $\pi_1$ and $\pi_2$ ramify over $2pd$ values, each of them marked with $p.$ Equivalently, the signature of the action of $H_j$ on $S$ is $(n_j; p, \stackrel{2dp}{\ldots}, p)$ where $n_j$ is the genus of $S/{H_j}$. Similarly as previously done,  the Riemann-Hurwitz formula ensures that $n_j=n$ for $j=1,2$. 

\s

In both cases, $H_1$ and $H_2$ are two $(p,n)$-gonal groups of $S,$ as desired.

\begin{rema1}
Note that if $d$ is odd then the $(p,n)$-gonal groups are conjugate,  but   if $d$ even then they are not.
\end{rema1}

\end{document}